\documentclass[final,leqno,onefignum,onetabnum]{siamltex}

\usepackage{amsmath}
\usepackage{amssymb}
\usepackage{amsbsy}
\usepackage{psfrag}
\usepackage[latin1]{inputenc}
\usepackage{graphicx}
\usepackage{cmmib57}

\renewcommand{\ldots}{\dotsc}

\newcommand{\real}{\mathbb{R}}
\newcommand{\diagonal}{\mathrm{diag}}

\newcommand{\figurewidth}{11.0cm}

\title{Explicit Time-Stepping for Stiff ODEs\thanks{Received by the editors
June 14, 2002; accepted for publication (in revised form) June~2, 2003; published
electronically December 5, 2003\@. \URL sisc/25-4/40962.html}}
\author{Kenneth Eriksson\thanks{Department of Computational Mathematics,
Chalmers University of Technology, SE--412 96 G\"oteborg, Sweden (kenneth@math.chalmers.se,
claes@math.chalmers.se, logg@math.chalmers.se).}
\and Claes Johnson\footnotemark[2]
\and Anders Logg\footnotemark[2]}

\begin{document}

\maketitle
\vspace{-1in}
\slugger{sisc}{2003}{25}{4}{1142--1157}
\vspace{.7in}

\setcounter{page}{1142}

\begin{abstract}
  We present a new strategy for solving stiff ODEs with explicit
  methods. By adaptively taking a small number of stabilizing small
  explicit time steps
  when necessary, a stiff ODE system can be stabilized enough to allow
  for time steps much larger than what is indicated by classical
  stability analysis. For many stiff problems the cost of the
  stabilizing small time steps is small, so the improvement
  is large.
  We illustrate the technique on a number of
  well-known stiff test problems.
\end{abstract}



\pagestyle{myheadings}
\thispagestyle{plain}
\markboth{KENNETH ERIKSSON, CLAES JOHNSSON, AND ANDERS LOGG}{EXPLICIT TIME-STEPPING FOR STIFF ODES}

\section{Introduction}

The classical wisdom developed in the 1950s regarding stiff ODEs
is that efficient
integration requires implicit (A-stable) methods, at
least outside of transients, where the time steps may be chosen
large from an accuracy point of view. Using an explicit method
(with a bounded stability region) the time steps have to
be small at all times for stability reasons,
and thus the advantage of a low cost per time step
may be counterbalanced by the necessity of taking a
large number of steps. As a result, the overall efficiency
of an explicit method for a stiff ODE may be small.

The question now is if it is possible to combine the low cost per time
step of an explicit method with the possibility of choosing large time
steps outside of transients. To reach this favorable combination, some
kind of stabilization of the explicit method seems to be needed, and
the basic question is then if the stabilization can be realized at a
low cost.

The stabilization technique proposed in this note relies on the
inherent property of the stiff problem itself, which is the rapid damping
of high frequencies. This is achieved by stabilizing the system
with a couple of small stabilizing (explicit Euler) steps.
We test this idea in adaptive form, where the size and
number of the small time steps are adaptively chosen according to the
size of different residuals. In particular, we compute rates of
divergence to determine the current mode $\lambda$ of largest
amplification and to determine a corresponding small time step $k\approx
\frac{1}{\lambda}$ with high damping. We test the
corresponding adaptive method in the setting of the $\mathrm{cG}(1)$
method
with fixed point iteration, effectively corresponding to an explicit
method if the number of iterations is kept small.
We show in a sequence of test problems that the
proposed adaptive method yields a significant reduction in work in comparison to
a standard implementation of an explicit method, with a typical speedup factor of $\sim 100$.
We conclude that, for many stiff problems, we may
efficiently use an explicit method, if only the explicit method is
adaptively stabilized with a relatively small number of small time steps,
and so we reach the desired combination of a low cost per time step
and the possibility of taking large time steps outside transients.
It remains to be seen if the proposed method can also be efficient in comparison to implicit methods.

We have been led to this question in our work on
multi-adaptive $\mathrm{cG}(q)$ or $\mathrm{dG}(q)$ ODE-solvers based on Galerkin's method with
continuous or discontinuous  polynomials of degree $q$, where individual
time steps are used for different components (see
\cite{logg:article:ma1,logg:article:ma2}).
These methods are implicit, and, to realize efficient implementations,
we need to use fixed point iteration with simple preconditioners.
With a limited (small) number of iterations, these iterative solvers
correspond to explicit time-stepping, and the same question of the cost
of stabilization arises.

The motivation for this work comes also from situations where for some reason
we are forced to using an explicit method, such as in simulations of
very large systems of chemical reactions or molecular dynamics, where
the formation of Jacobians and the solution of the linear system become very
expensive.

Possibly, similar techniques may be used also to stabilize
multigrid smoothers.

\section{Basic strategy}
\label{sec:testequation}

We consider first the {\em test equation\/}: Find $u:[0,\infty
)\rightarrow \real$ such that
\begin{equation}\label{testequation}
\begin{split}
\dot u(t)+\lambda u(t)&=0\quad\mbox{for }t>0,\\
u(0)&=u^0,
\end{split}
\end{equation}
where $\lambda> 0$ and $u^0$ is a given initial condition. The solution is given by $u(t)=\exp(-\lambda t)u^0$.
We define the {\em transient\/} as $\{t>0:\lambda t\le C\}$ with $C$ a moderate constant. Outside the
transient, $u(t)=\exp(-\lambda t)u^0$ will be small.

The {\em explicit Euler method\/} for the test equation reads
\[
U^n=-k_n\lambda U^{n-1}+U^{n-1}=(1-k_n\lambda )U^{n-1}.
\]
This method is conditionally stable and requires that $k_n\lambda\le 2$, which, outside of transients,
is too restrictive for large $\lambda$.

Now let $K$ be a large time step satisfying $K\lambda >2$ and let $k$ be a small time step chosen so that
$k\lambda <2$. Consider the method
\begin{equation}\label{explimplE0}
U^n=(1-k\lambda )^{m}(1-K\lambda )U^{n-1},
\end{equation}
corresponding to one explicit Euler step with large time step $K$ and $m$
explicit Euler steps with small time steps $k$, where $m$ is a positive
integer to be determined. Altogether this corresponds to a time step
of size $k_n=K+mk$. Defining the polynomial function $P(x)=(1-\theta
x)^{m}(1-x)$, where $\theta =\frac{k}{K}$, we can write method
(\ref{explimplE0}) in the form
\[
U^n=P(K\lambda )U^{n-1}.
\]
We now seek to choose $m$ so that $\vert P(K\lambda )\vert\le 1$,
which is needed for stability. We thus need to satisfy
\[
|1-k\lambda|^m (K\lambda-1)\le 1,
\]
that is,
\begin{equation}
\label{eq:m}
m\ge \frac{\log(K\lambda -1)}{-\log|1-k\lambda|}\approx \frac{\log(K\lambda )}{c},
\end{equation}
with $c=k\lambda$ a moderate constant; for definiteness, let $c=1/2$.

We conclude that  $m$ will be of moderate size, and consequently only a small fraction
of the total time interval will be spent on time-stepping with the small time step $k$.
To see this, define the \emph{cost} as
\[
\alpha = \frac{1+m}{K+km} \in (1/K,1/k),
\]
that is, the number of time steps per unit time interval.
Classical stability analysis with $m=0$ gives
\begin{equation}
	\label{eq:alpha0}
	\alpha = 1/k_n = \lambda/2,
\end{equation}
with a maximum time step $k_n=K=2/\lambda$.
Using (\ref{eq:m}) we instead find
\begin{equation}
	\alpha \approx \frac{1 + \log(K\lambda)/c}{K + \log(K\lambda)/\lambda} \approx
	\frac{\lambda}{c} \log(K\lambda)/(K\lambda) \ll \lambda/c
\end{equation}
for $K\lambda \gg 1$.
The cost is thus decreased by the cost reduction factor
\[
	\frac{2\log(K\lambda)}{cK\lambda}
\sim \frac{\log(K\lambda)}{K\lambda},
\]
which can be quite significant for large values of $K\lambda$.

A similar analysis applies to the system $\dot{u} + Au = 0$ if the eigenvalues $\{\lambda_i\}_{i=1}^N$ of
$A$ are separated with $0 \leq \lambda_1 \leq \cdots \leq \lambda_{i-1} \leq 2/K$ and
$2/K \ll \lambda_i \leq \cdots \leq \lambda_N$. In this case, the cost is decreased by
a factor of
\begin{equation}
	\label{eq:gain,li}
	\frac{2\log(K\lambda_i)}{c K\lambda_i} \sim \frac{\log(K\lambda_i)}{K\lambda_i}.
\end{equation}

In recent independent work by Gear and Kevrekidis~\cite{GeaKev02}, a
similar idea of combining small and large time steps for a class of
stiff problems with a clear separation of slow and fast time scales
is developed. That work, however, is not focused on adaptivity to the same
extent as ours, which does not require any a~priori information about
the nature of the stiffness (for example, the distribution of eigenvalues).

\section{Parabolic problems}
\label{sec:parabolic}

For a parabolic system,
\begin{equation}\label{eq:parabolic}
\begin{split}
\dot u(t)+A u(t)&=0\quad\mbox{for }t>0,\\
u(0)&=u^0,
\end{split}
\end{equation}
with $A$ a symmetric positive semidefinite $N\times N$ matrix
and $u^0\in \real^N$ a given initial condition,
the basic damping strategy outlined in the previous section may fail to be
efficient. This can be seen directly from (\ref{eq:gain,li});
for efficiency we need $K\lambda_i \gg 2$, but with
the eigenvalues of $A$ distributed over the interval $[0,\lambda_N]$,
for example with $\lambda_i \sim i^2$ as for a typical (one-dimensional) parabolic problem,
one cannot have both $\lambda_{i-1} \leq 2/K$ and $K\lambda_i \gg 2$!

We conclude that, for a parabolic problem, the sequence of time steps,
or equivalently the polynomial $P(x)$, has to be chosen differently.
We thus seek a more general sequence $k_1,\ldots,k_m$ of time steps such that
$\vert P(x) \vert \leq 1$ for $x\in [0, K\lambda_N]$, with
\begin{displaymath}
P(x) = \left( 1 - \frac{k_1 x}{K} \right) \cdots \left( 1 - \frac{k_m x}{K} \right)
\end{displaymath}
and $K$ a given maximum step size.

\subsection{Chebyshev damping}

A good candidate for $P$ is given by the shifted Chebyshev polynomial of
degree $m$,
\begin{displaymath}
  P_c(x) = T_m \left( 1 - \frac{2x}{K\lambda_N} \right).
\end{displaymath}
This gives $P_c(0) = 1$ and $\vert P_c(x) \vert \leq 1$ on $[0,K\lambda_N]$ (see
Figure~\ref{fig:comparison}).
A similar approach is taken in~\cite{Ver96}.

Analyzing the cost as before, we have $\alpha = m / (k_1 + \cdots + k_m)$,
with
\begin{displaymath}
 	k_i = 2/(\lambda_N (1-s_{m+1-i})),
\end{displaymath}
and $s_i$ the $i$th zero of the Chebyshev polynomial $T_m=T_m(s)$,
given by $s_i = \cos((2i-1)\pi/(2m))$, $i=1,\ldots,m$. It follows that
\begin{displaymath}
	\alpha = \frac{m \lambda_N/2}{1/(1-s_1) + \cdots + 1/(1-s_m)} =
   \frac{m \lambda_N/2}{m^2} = \lambda_N/2m.
\end{displaymath}
The reduction in cost compared to $\lambda_N/2$ is thus a factor $1/m$.
A restriction on the maximum value of $m$ is given by $k_m = 2/(\lambda_N (1-s_1)) \leq K$;
that is,
\begin{displaymath}
K \geq \frac{2}{\lambda_N (1-\cos(\pi/(2m)))} \approx \frac{2}{\lambda_N \pi^2/(8m^2)} = 16m^2 / (\lambda_N \pi^2)
\end{displaymath}
for $m \gg 1$. With $m = (\pi/4) \sqrt{K\lambda_N}$, the cost reduction factor for Chebyshev damping
is thus given by
\begin{displaymath}
	1/m = \frac{4}{\pi \sqrt{K\lambda_N}} \sim (K\lambda_N)^{-1/2}.
\end{displaymath}

\subsection{Dyadic damping}

Another approach is a modified version of the simple approach of small and large time steps discussed
in the previous section. As discussed above, the problem with this basic method is that the damping is
inefficient for intermediate eigenmodes. To solve this problem, we modify
the simple method of large and small time steps to include also time steps of intermediate size.
A related approach is taken in~\cite{GeaKev03}.
Starting with a couple of small time steps of size $k$, we thus increase the time step
gradually by (say) a factor of~$2$, until we reach the largest step size $K$.
Correspondingly, we decrease the number of time steps at each level by a factor of~2.
The polynomial $P$ is now given by
\begin{displaymath}
  P_d(x) = \prod_{j=q+1}^p \left(1-\frac{2^j x}{K\lambda_N}\right)
  \prod_{i=0}^q \left(1-\frac{2^i x}{K\lambda_N} \right)^{2^{q-i}},
\end{displaymath}
where $p$ and $q\leq p$ are two integers to be determined.
We thus obtain a sequence of increasing time steps, starting with $2^q$ time steps of size
$k = 1/\lambda_N$, $2^{q-1}$ time steps of size $2 k$, and so on, until we reach
level $q$ where we take one time step of size $2^q k$. We then continue to
increase the size of the time step, with only one time step at each level, until we reach
a time step of size $2^p k = K$.

With $p$ given, we determine the minimal value of $q$ for $|P_d(x)| \leq 1$
on $[0,K\lambda_N]$ for $p=0,1,\ldots$, and find $q(p) \approx \frac{2}{3}(p-1)$;
see Table~\ref{tab:qp}. The cost is now given by
\begin{eqnarray*}
  \alpha &=& \frac{(2^q + \cdots + 1) + (p - q)}{\frac{1}{\lambda_N}(2^q(q+1) + (2^{q+1} + \cdots + 2^p))}
  = \frac{\lambda_N ( (2^{q+1} - 1) + (p - q))}{2^q(q+1) + (2^{p+1}-2^{q+1})} \\
  &\approx& \frac{\lambda_N/2}{2^{p-q-1}} \approx
  \frac{\lambda_N}{2} 2^{-p+2(p-1)/3+1} = \frac{\lambda_N}{2} 2^{-(p-1)/3}.
\end{eqnarray*}
Since $p = \log(K\lambda_N)/\log(2)$, the reduction in cost for dyadic damping is then a factor
\begin{displaymath}
\frac{1}{(K\lambda_N/2)^{1/3}} \sim (K\lambda_N)^{-1/3},
\end{displaymath}
which is competitive with the optimal cost reduction of Chebyshev damping.

\begin{table}
\footnotesize
\caption{Number of required levels with multiple zeros, $q$, as a function of the total number of levels, $p$.}\label{tab:qp}
\begin{center}
\begin{tabular}{|l|rrrrrrrrrrrrrrrrr|}
      \hline
      $p$ & 0 & 1 & 2 & 3 & 4 &  5 & 6 & 7 & 8 & 9 & 10 & 11 & 12 & 13 & 14 & 15 & 16 \\
      \hline
      $q$ & 0 & 0 & 0 & 1 & 2 &  3 & 3 & 4 & 4 & 5 &  6 &  7 &  8 &  8 &  9 & 10 & 10 \\
      \hline
\end{tabular}
\end{center}
\end{table}

\section{Comparison of methods}

We summarize our findings so far as follows. The reduction in cost for the simple approach of large
unstable time steps followed by small stabilizing time steps (outlined in section~\ref{sec:testequation})
is a factor $\log(K\lambda_N)/(K\lambda_N)$, and thus this simple approach can be much more efficient than
both of the two approaches, Chebyshev damping and dyadic damping, discussed in the previous section.
This however requires that the problem in question has a gap in its eigenvalue spectrum; that is, we must have
a clear separation into small (stable) eigenvalues and large (unstable) eigenvalues.

In the case of a parabolic problem, without a gap in the eigenvalue spectrum, the simple approach of
section~\ref{sec:testequation} will not work. In this case, Chebyshev damping and dyadic damping give a reduction
in cost which is a factor $(K\lambda_N)^{-1/2}$ or $(K\lambda_N)^{-1/3}$, respectively.
The efficiency for a parabolic problem is thus comparable for the two methods. Since the
dyadic damping method is a slightly modified version of the simple approach of section~\ref{sec:testequation},
consisting of gradually and carefully increasing the time step size after stabilization with the smallest
time step, thus being simple and flexible, we believe this method to be advantageous over the Chebyshev method.

\begin{figure}[ht]
\begin{center}
    \psfrag{x}{$x$}
    \psfrag{PC}{$P_d(x)$}
    \psfrag{PD}{$P_c(x)$}
    \includegraphics[width=25pc]{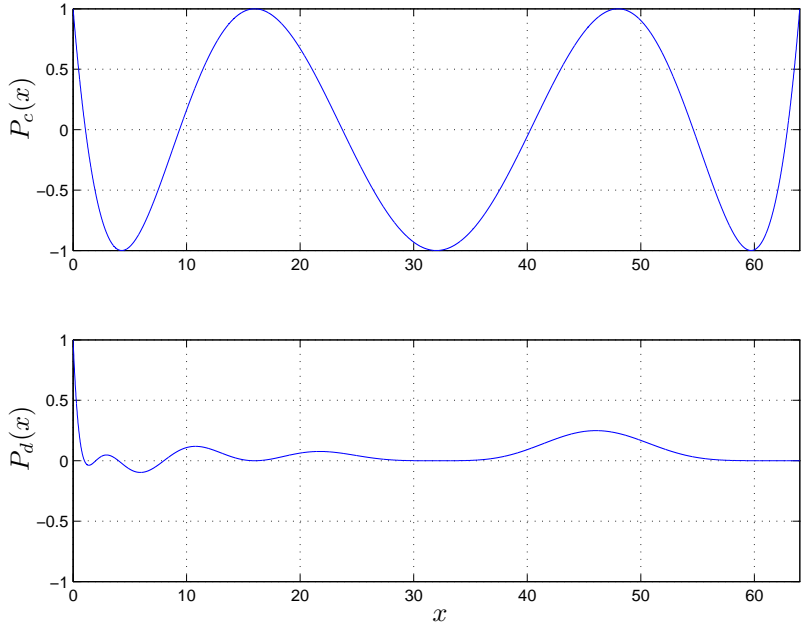}
\end{center}
    \caption{A comparison of the two polynomials for $K\lambda_N=64$:
      we take $m=6$ for the shifted Chebyshev polynomial $P_c(x)$ and\/
      $(p,q)=(6,3)$ for the dyadic polynomial $P_d(x)$.
      With $K=1$, the costs are\/ $5.3$ and\/~$8$, respectively.}\label{fig:comparison}
\end{figure}

As a comparison, we plot in Figure~\ref{fig:comparison} the shifted Chebyshev polynomial
$P_c(x)$ and the polynomial $P_d(x)$ of dyadic damping for $K\lambda_N = 64$.
With $K\lambda_N = 64$, we obtain $m = (\pi/4)\sqrt{K\lambda_N} \approx 6$ for the Chebyshev polynomial.
For the dyadic damping polynomial, we have $p=\log(K\lambda_N)/\log(2) = 6$, and from Table~\ref{tab:qp}, we obtain $q=3$.

\begin{figure}[ht]
\begin{center}
    \psfrag{x}{$x$}
    \psfrag{PC}{$P_d(x)$}
    \psfrag{PD}{$P_c(x)$}
    \includegraphics[width=20pc]{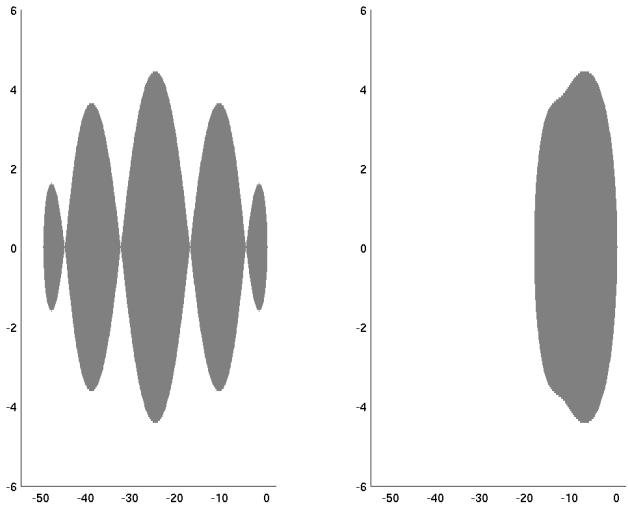}
\end{center}
    \caption{A comparison of stability regions for the two polynomials,
      with $m=5$ for the shifted Chebyshev polynomial $P_c(z)$ (left) and\/ $(p,q)=(3,1)$ for the dyadic polynomial
      $P_d(z)$ (right).}\label{fig:stabilityregions}
\end{figure}

As another comparison, we plot
in Figure~\ref{fig:stabilityregions} the stability regions for the two polynomials $P_c(z)$ and $P_d(z)$
for $z\in\mathbb{C}$.
In this context, the two polynomials are given by
\begin{equation}
  P_c(z) = \prod_{i=1}^m \left(1+\frac{z/m^2}{1-s_i} \right)
\end{equation}
for $s_i = \cos((2i-1)\pi/(2m))$ and
\begin{equation}
  P_d(z) = \prod_{j=q+1}^p \left(1+\theta_j z\right)
           \prod_{i=0}^q   \left(1+\theta_i z\right)^{2^{q-1}}
\end{equation}
for $\theta_i = 2^i / (2^q(q+1) + 2^{p+1} - 2^{q+1})$.
Note that we take $z = - \bar{K} \lambda$, where $\bar{K}=k_1+\cdots + k_m$ is the sum of all time steps in the sequence,
rather than $z = - K \lambda$, where $K \leq \bar{K}$ is the maximum time step in the sequence.
To make the comparison fair, we
take $m=5$ for the Chebyshev polynomial and $(p,q)=(3,1)$ for the dyadic damping polynomial, which gives
the same number of time steps ($m=2^{q+1}-1+p-q=5$) for the two methods.

From Figures \ref{fig:comparison} and~\ref{fig:stabilityregions}, we see that Chebyshev damping can be slightly
more efficient than dyadic damping; the cost is smaller for a given value of $K\lambda_N$ (Figure~\ref{fig:comparison}) and
the stability interval on the negative real axis is larger (Figure~\ref{fig:stabilityregions}). However, with dyadic damping we don't need to
assume that the eigenvalues of $A$ in (\ref{eq:parabolic}) lie in a narrow
strip along the negative real axis in the complex plane, as is needed with
Chebyshev damping. (The stability region of the Chebyshev polynomial can be slightly modified to
include a narrow strip along the negative real axis; see~\cite{Ver96}.)

\section{The general nonlinear problem}
\label{sec:nonlinear}

We consider now the general nonlinear problem,
\begin{equation}
	\label{eq:u'=f}
	\begin{split}
		\dot{u}(t) &= f(u(t)) \quad \mbox{for } t>0,\\
		u(0) &= u^0,
	\end{split}
\end{equation}
where $f:\real^N\times(0,T]\rightarrow\real^N$ is
a given bounded and differentiable function.
The explicit Euler method for (\ref{eq:u'=f}) reads
\begin{equation} \label{eq:U}
	U^n = U^{n-1} + k_n f(U^{n-1}),
\end{equation}
where the function $U=U(t)$ is piecewise constant and right-continuous with
$U^n = U(t_n^+)$.
The exact solution $u$
satisfies a similar relation,
\begin{equation} \label{eq:u}
	u^n = u^{n-1} + \int_{t_{n-1}}^{t_n} f(u(t))\ dt,
\end{equation}
with $u^n = u(t_n)$.
Subtracting (\ref{eq:U}) and (\ref{eq:u}), we find that the error $e(t) = U(t) - u(t)$ satisfies
\begin{displaymath}
	e(t_n^+) - e(t_{n-1}^+)
	 =  \int_{t_{n-1}}^{t_n} J (U(t) - u(t)) \ dt
	 =  \int_{t_{n-1}}^{t_n} J e \ dt,
\end{displaymath}
where $J$ is the Jacobian $\frac{\partial f}{\partial u}$ of the right-hand side evaluated
at a mean value of $U$ and $u$.
We conclude that the efficiency of the proposed method for the general nonlinear problem
will be determined by the distribution of the eigenvalues of the Jacobian.

\section{Iterative methods}
\label{sec:iterative}

From another viewpoint, we may consider using an explicit-type iterative
method for solving the discrete equations arising from an implicit method.
The implicit $\mathrm{cG}(1)$ method with midpoint quadrature for the general nonlinear problem
(\ref{eq:u'=f}) reads
\begin{equation}
	\label{eq:cg1}
	U^n = U^{n-1} + k_n f\left(\frac{U^{n-1}+U^n}{2}\right).
\end{equation}
We can solve this system of nonlinear equations for $U^n$ using Newton's
method, but the simplest and cheapest method is to apply fixed point
iteration directly to (\ref{eq:cg1}); that is,
\begin{displaymath}
	U^{nl} = U^{n-1} + k_n f\left(\frac{U^{n-1}+U^{n,l-1}}{2}\right),
\end{displaymath}
for $l=1,2,\ldots$ until convergence occurs with $U^{n,0} = U^{n-1}$. The
fixed point
iteration converges for a small enough time step $k_n$, so
the stability condition for a standard explicit method appears also in explicit-type
iterative methods as a condition for convergence of the iterative solution of the implicit
equations.

To determine a stop criterion for the fixed point iteration, we measure the size of
the \emph{discrete residual},
\begin{displaymath}
	r^{nl} = \frac{1}{k_n}(U^{nl} - U^{n-1}) - f\left(\frac{U^{n-1}+U^{nl}}{2}\right),
\end{displaymath}
which should be zero for the true $\mathrm{cG}(1)$ approximation. We continue
the iterations until the discrete residual is smaller than some tolerance $\mathrm{tol}>0$.
Usually only a couple of iterations are needed.
Estimating the error $e(T)=U(T)-u(T)$ at final time $T$
(see~\cite{logg:article:ma1}), we have
\[
	\|e(T)\| \leq S(T) \max_{[0,T]} k \|R\| + S^0(T) \max_{[0,T]} \|r\|,
\]
where $R(t)=\dot{U}(t) - f(U(t))$ is the \emph{continuous residual} and
$S(T)$ and $S^0(T)$ are stability factors. For the test equation, we have
$S(T) \leq 1$ and $S^0(T) \leq 1/\lambda$, which suggests that for a
typical stiff problem we can take $\mathrm{tol} = \mathrm{TOL}$, where
$\mathrm{TOL}$ is a tolerance for the error $e(T)$ at final time.

For the discrete residual, we have
\begin{eqnarray*}
		r^{nl}
		&=& \frac{1}{k_n}(U^{nl} - U^{n-1}) - f\left(\frac{U^{n-1}+U^{nl}}{2}\right) \\
		&=& f\left(\frac{U^{n-1}+U^{n,l-1}}{2}\right) - f\left(\frac{U^{n-1}+U^{nl}}{2}\right)
		 =  J \frac{U^{n,l-1}-U^{nl}}{2} \\
		&=& \frac{1}{2}J\left[ U^{n,l-1}-U^{n-1}-k_nf\left(\frac{U^{n-1}+U^{n,l-1}}{2}\right) \right],
\end{eqnarray*}
which gives
\begin{equation}
	\label{eq:convergence}
	r^{nl} = \frac{k_n}{2} J r^{n,l-1}.
\end{equation}

Thus, by measuring the size of the discrete residual, we obtain
information about the stiff nature of the problem, in particular the
eigenvalue of the current unstable eigenmode, which can be used in an
adaptive algorithm targeted precisely at stabilizing the current
unstable eigenmode. We discuss this in more detail below in section~\ref{sec:algorithm}.

\section{Multi-adaptive solvers}

In a multi-adaptive solver, we use individual time steps for different components.
An important part of the algorithm described
in \cite{logg:article:ma1,logg:article:ma2}
is the iterative fixed point solution of the discrete equations on time slabs.
The simple strategy to take small damping steps to stabilize the system applies
also in the multi-adaptive setting, where we may also target individual eigenmodes
(if these are represented by different components) using individual damping steps.
This will be explored further in another note.

\section{An adaptive algorithm}
\label{sec:algorithm}

The question is now whether we can choose the time step sequence automatically in an
adaptive algorithm. We approach this question in the setting of an implicit
method combined with an explicit-type iterative solver as in section~\ref{sec:iterative}. We give the algorithm in the case of the simple damping strategy
outlined in section~\ref{sec:testequation}, with an extension to
parabolic problems as described in section~\ref{sec:parabolic}.

A simple adaptive algorithm for the standard $\mathrm{cG}(1)$ method with
iterative solution of the discrete equations reads as follows.
\begin{enumerate}
\item
	Determine a suitable initial time step $k_1$, and solve the discrete
	equations for the solution $U(t)$ on $(t_0,t_1)$.
\item
	Repeat on $(t_{n-1},t_n)$ for $n=2,3,\ldots$ until $t_n \geq T$:
	\begin{enumerate}
		\item
		Evaluate the continuous residual $R_{n-1}$ from the previous time interval.
		\item
		Determine the new time step $k_n$ based on $R_{n-1}$.
		\item
		Solve the discrete equations on $(t_{n-1},t_n)$ using fixed point iteration.
	\end{enumerate}
\end{enumerate}
In reality we want to control the global error, which means we also have to
solve the dual problem, compute stability factors (or weights), evaluate
an a~posteriori error estimate, and possibly repeat the process until the error
is below a given tolerance $\mathrm{TOL}>0$. The full algorithm is thus slightly
more elaborate, but the basic algorithm presented here is the central part.
See~\cite{EriEst95} for a discussion.

We comment also on step~2(b): For the $\mathrm{cG}(1)$ method we would like
to take $k_n = \mathrm{TOL} / (S(T) \|R_{n-1}\|)$, but this introduces oscillations
in the size of the time step. A small residual gives a large time step which results
in a large residual, and so on. To avoid this, the simple step size selection has to
be combined with a regulator of some kind; see \cite{GusLun88,Sod98} or~\cite{logg:article:ma2}. It turns out that a simple strategy that works well
in many situations is to take $k_n$ as the geometric mean value
\begin{equation}
	\label{eq:regulator}
	k_n = \frac{2 \tilde{k}_n k_{n-1}}{\tilde{k}_n + k_{n-1}},
\end{equation}
where $\tilde{k}_n = \mathrm{TOL} / (S(T) \|R_{n-1}\|)$.

Now, for a stiff problem, what may go wrong is step~2(c); if the time step $k_n$ is
too large, the fixed point iterations will not converge. To be able to handle stiff
problems using the technique discussed above, we propose the following simple
modification of the adaptive algorithm.

\begin{enumerate}
\item
	Determine a suitable initial time step $k_1$, and solve the discrete
	equations for the solution $U(t)$ on $(t_0,t_1)$.
\item
	Repeat on $(t_{n-1},t_n)$ for $n=2,3,\ldots$ until $t_n \geq T$:
	\begin{enumerate}
		\item
		Evaluate the continuous residual $R_{n-1}$ from the previous time interval.
		\item
		Determine the new time step $k_n$ based on $R_{n-1}$.
		\item
		Solve the discrete equations on $(t_{n-1},t_n)$ using fixed point iteration.
		\item
		If 2(c) did not work, compute
		\begin{displaymath}
			L = \frac{2}{k_n}\frac{\|r^l\|}{\|r^{l-1}\|},
		\end{displaymath}
		and take $m=\log(k_n L)$
		explicit Euler steps with time step $k=c/L$ and $c$ close to $1$.
		\item
		Try again starting at 2(a) with $n\rightarrow n+m$.
	\end{enumerate}
\end{enumerate}

In the analysis of section~\ref{sec:testequation},
we had $c=1/2$, but it is clear that the damping steps
will be more efficient if we have $c$ close to~$1$.
An implementation of this algorithm in the form of a simple MATLAB code
is available for inspection~\cite{logg:www:software}, including among others
the test problems presented in the next section.

We also note that, by (\ref{eq:regulator}), we have $k_n \leq 2 k_{n-1}$, so,
following the sequence of small stabilizing time steps, the size of the time step
will be increased gradually, doubling the time step until we reach $k_n \sim K$
or the system becomes unstable again, whichever comes first. This automatically gives a sequence
of time steps similar to that of dyadic damping described in section~\ref{sec:parabolic},
with the difference being that most of the damping is made with the smallest time step.
For a parabolic problem, we modify the algorithm to increase the time step more carefully,
as determined by Table~\ref{tab:qp}, with $p$ given by $p=\log(k_n L)/\log(2)$.

\section{Examples}

To illustrate the technique, we take a simple standard implementation of the
$\mathrm{cG}(1)$-method (with an explicit fixed point solution of the discrete equations)
and add a couple of lines to handle the stabilization of stiff problems.
We try this code on a number of well-known stiff problems taken from the ODE literature
and conclude that we are able to handle stiff problems with this explicit code.

When referring to the cost $\alpha$ below, this also includes the
number of fixed point
iterations needed to compute the $\mathrm{cG}(1)$ solution on intervals where the
iterations converge. This is compared to the cost $\alpha_0$ for the standard
$\mathrm{cG}(1)$ method in which we are forced to take small time
steps all the time. (These small time steps are marked by dashed lines.)
For all example problems below, we report both the cost
$\alpha$ and the cost reduction factor $\alpha/\alpha_0$.

\subsection{The test equation}

\begin{figure}[t]
	\psfrag{t}{$t$}
	\psfrag{U}{$U(t)$}
	\psfrag{k}{$k(t)$}
\begin{center}
\includegraphics[width=\figurewidth]{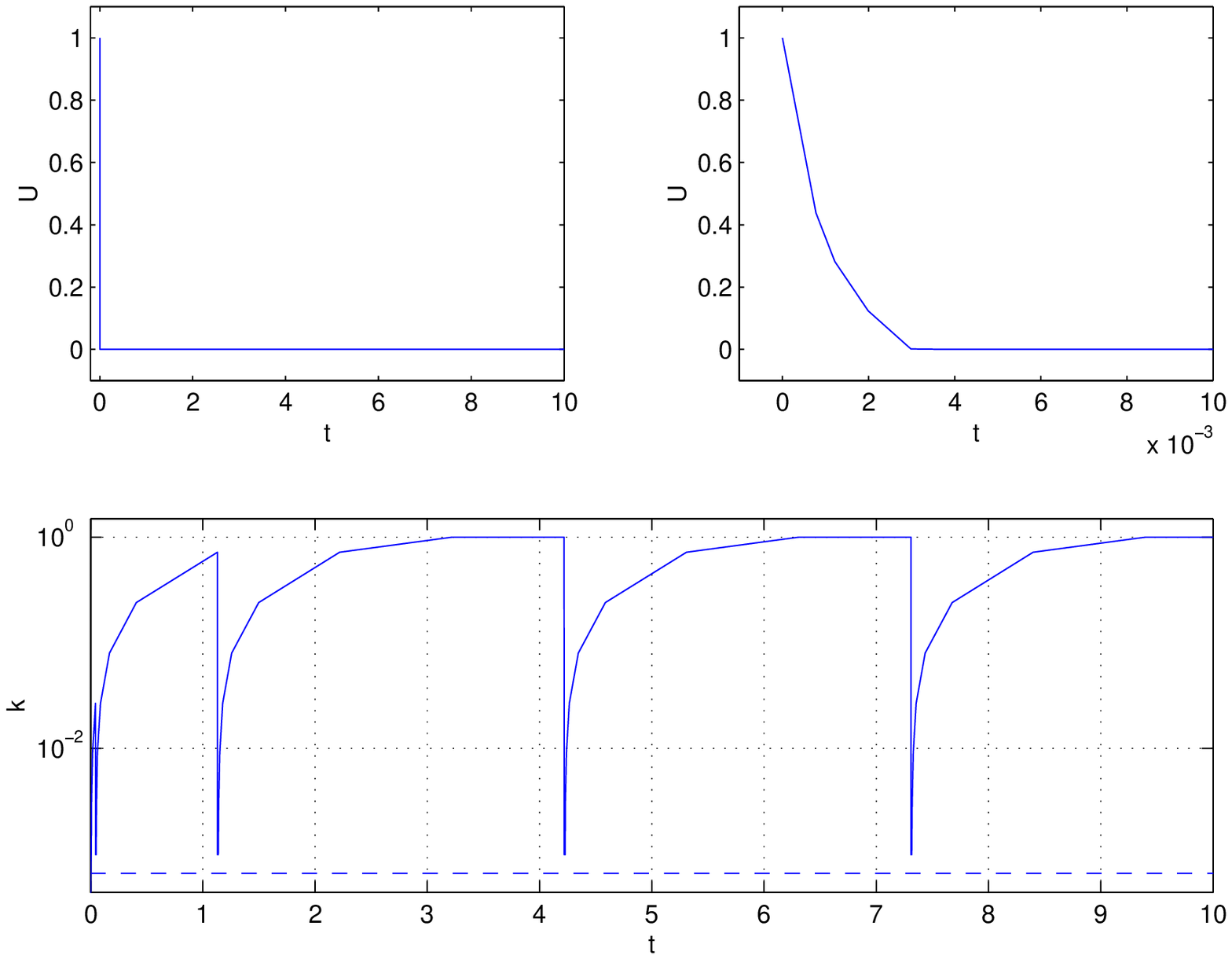}
\end{center}
\caption{Solution and time step sequence for\/ {\rm(\ref{eq:problem1})},
$\alpha/\alpha_0 \approx 1/310$.}\label{fig:problem1}
	\end{figure}

The first problem we try is the test equation,
\begin{equation}
	\label{eq:problem1}
	\begin{split}
		\dot u(t)+\lambda u(t)&=0\quad\mbox{for }t>0,\\
		u(0)&=u^0,
	\end{split}
\end{equation}
on $[0,10]$, where we choose $u^0=1$ and $\lambda=1000$. As is evident from
Figure \ref{fig:problem1}, the time step sequence is automatically chosen
in agreement with the previous discussion. The cost is only
$\alpha\approx 6$
and the cost reduction factor is $\alpha/\alpha_0\approx 1/310$.

Note how the time steps are drastically decreased (step~2(d) in the adaptive algorithm)
when the system needs to be stabilized and then gradually increased until stabilization again
is needed.

\subsection{The test system}

\begin{figure}[t]
	\psfrag{t}{$t$}
	\psfrag{U}{$U(t)$}
	\psfrag{k}{$k(t)$}
\begin{center}
\includegraphics[width=\figurewidth]{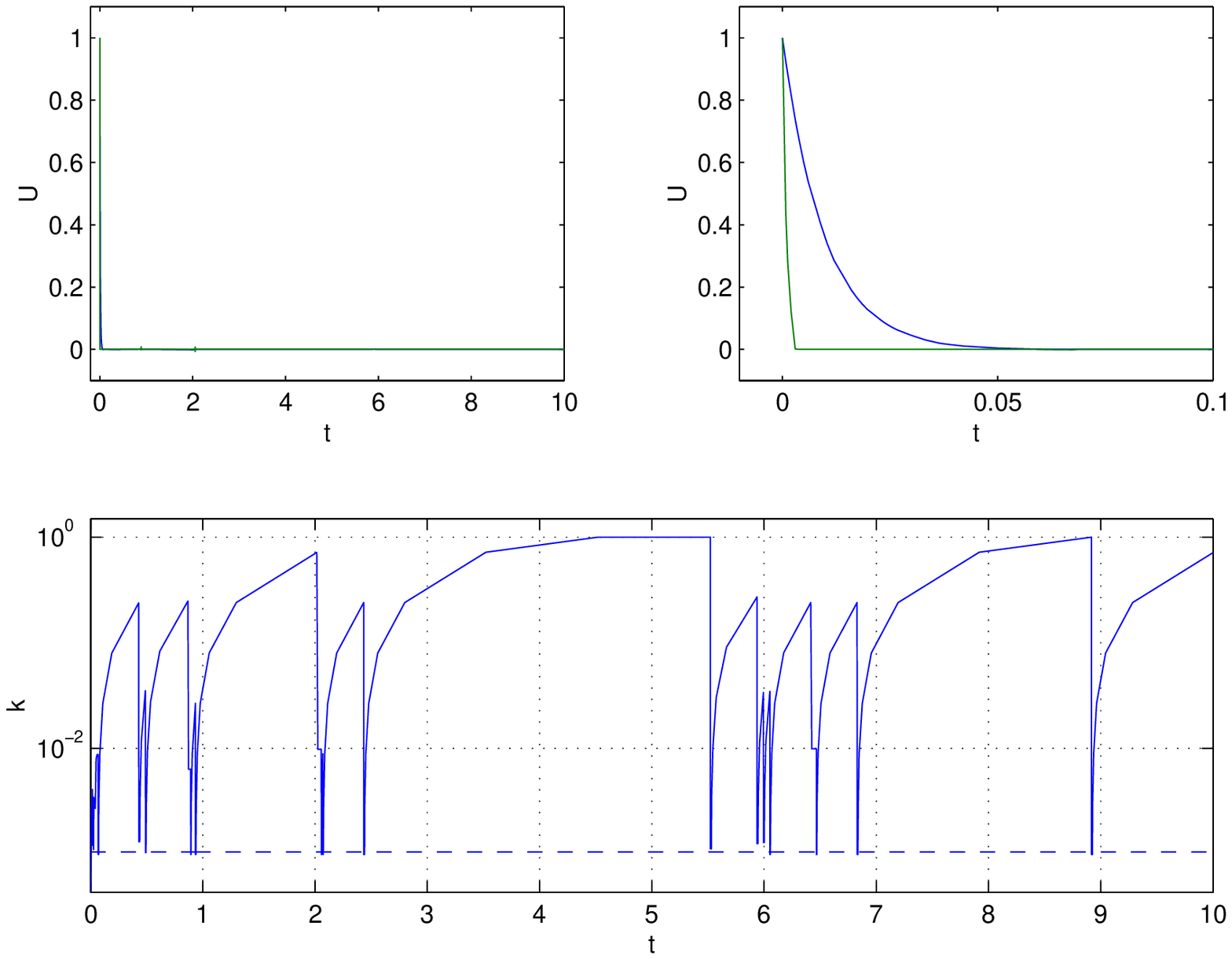}
\end{center}
\caption{Solution and time step sequence for\/ {\rm(\ref{eq:problem2})},
$\alpha/\alpha_0 \approx 1/104$.}\label{fig:problem2}
\end{figure}

For the test system,
\begin{equation}
	\label{eq:problem2}
	\begin{split}
		\dot u(t)+Au(t)&=0\quad\mbox{for }t>0,\\
		u(0)&=u^0,
	\end{split}
\end{equation}
on $[0,10]$, we take $A = \diagonal(100,1000)$ and $u^0=(1,1)$.
There are now two eigenmodes with large
eigenvalues that need to be damped out.
The dominant eigenvalue of $A$ is $\lambda_2=1000$, and most
of the damping steps are chosen to damp out this eigenmode, but some of the damping
steps are chosen based on the second largest eigenvalue $\lambda_1=100$ (see Figure
\ref{fig:problem2}). When to damp out
which eigenmode is automatically decided by the adaptive algorithm; the bad eigenmode that
needs to be damped out becomes visible in the iterative solution
process.
Since there is an additional eigenvalue,
the cost is somewhat larger than for the scalar test problem,
$\alpha\approx 18$, which gives a
cost reduction factor of $\alpha/\alpha_0\approx 1/104$.

\subsection{A linear nonnormal problem}

\begin{figure}[t]
	\psfrag{t}{$t$}
	\psfrag{U}{$U(t)$}
	\psfrag{k}{$k(t)$}
\begin{center}
\includegraphics[width=\figurewidth]{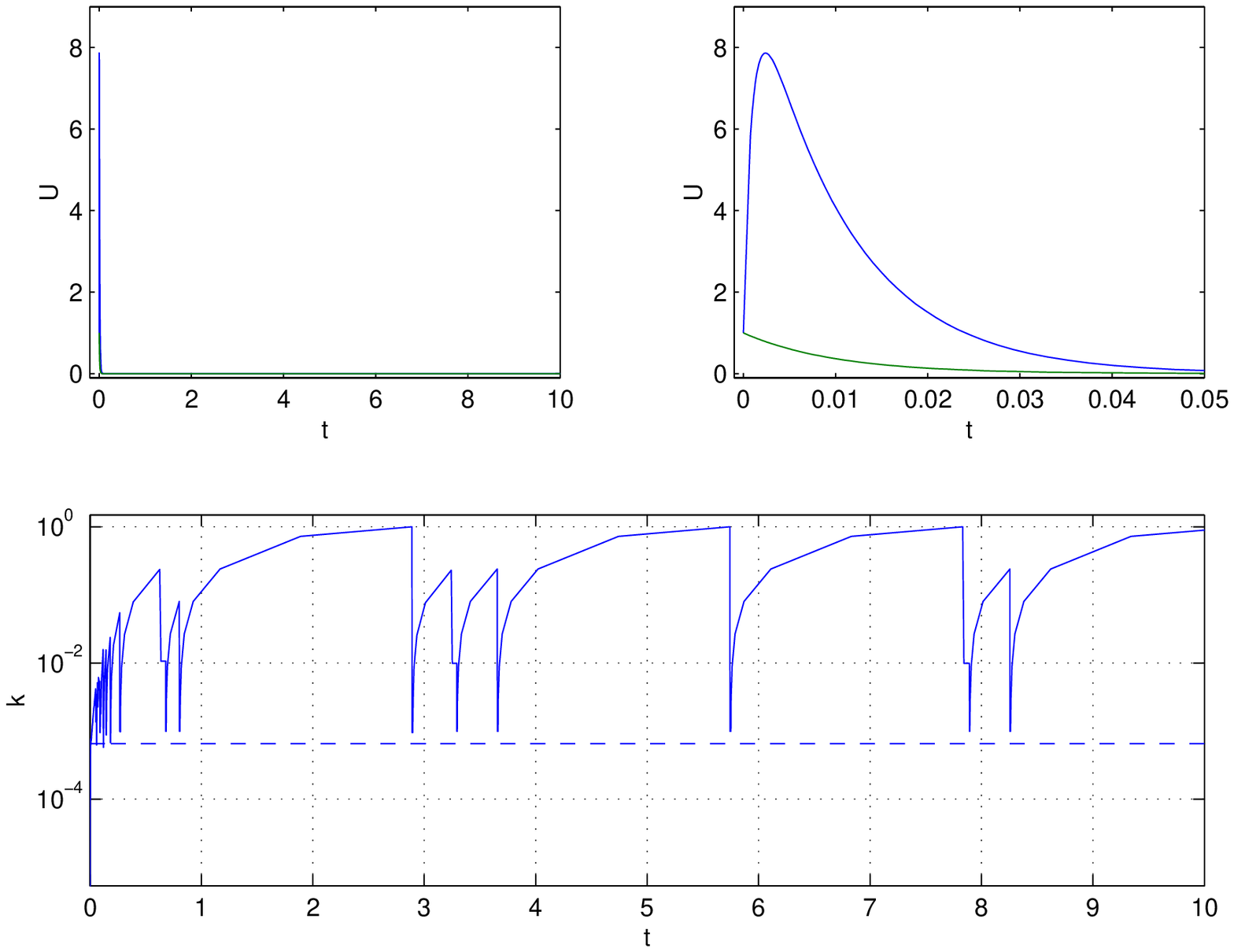}
\end{center}
\caption{Solution and time step sequence for\/ {\rm(\ref{eq:problem4})},
$\alpha/\alpha_0 \approx 1/180$.}\label{fig:problem4}
\end{figure}

The method behaves similarly even if we
make the matrix $A$ highly nonnormal. We now solve
\begin{equation}
	\label{eq:problem4}
	\begin{split}
		\dot u(t)+Au(t)&=0\quad\mbox{for }t>0,\\
		u(0)&=u^0,
	\end{split}
\end{equation}
on $[0,10]$, with
\begin{displaymath}
	A =
	\left[
	\begin{array}{cc}
		 1000 & -10000 \\
       0    & 100
	\end{array}
	\right]
\end{displaymath}
and $u^0=(1,1)$.
The cost is about the same as for the previous problem,
$\alpha\approx 17$, but the cost reduction factor is better:
$\alpha/\alpha_0\approx 1/180$. The solution and the time step sequence are given in Figure
\ref{fig:problem4}.

\subsection{The HIRES problem}

\begin{figure}[ht]
	\psfrag{t}{$t$}
	\psfrag{U}{$U(t)$}
	\psfrag{k}{$k(t)$}
\begin{center}
\includegraphics[width=24pc]{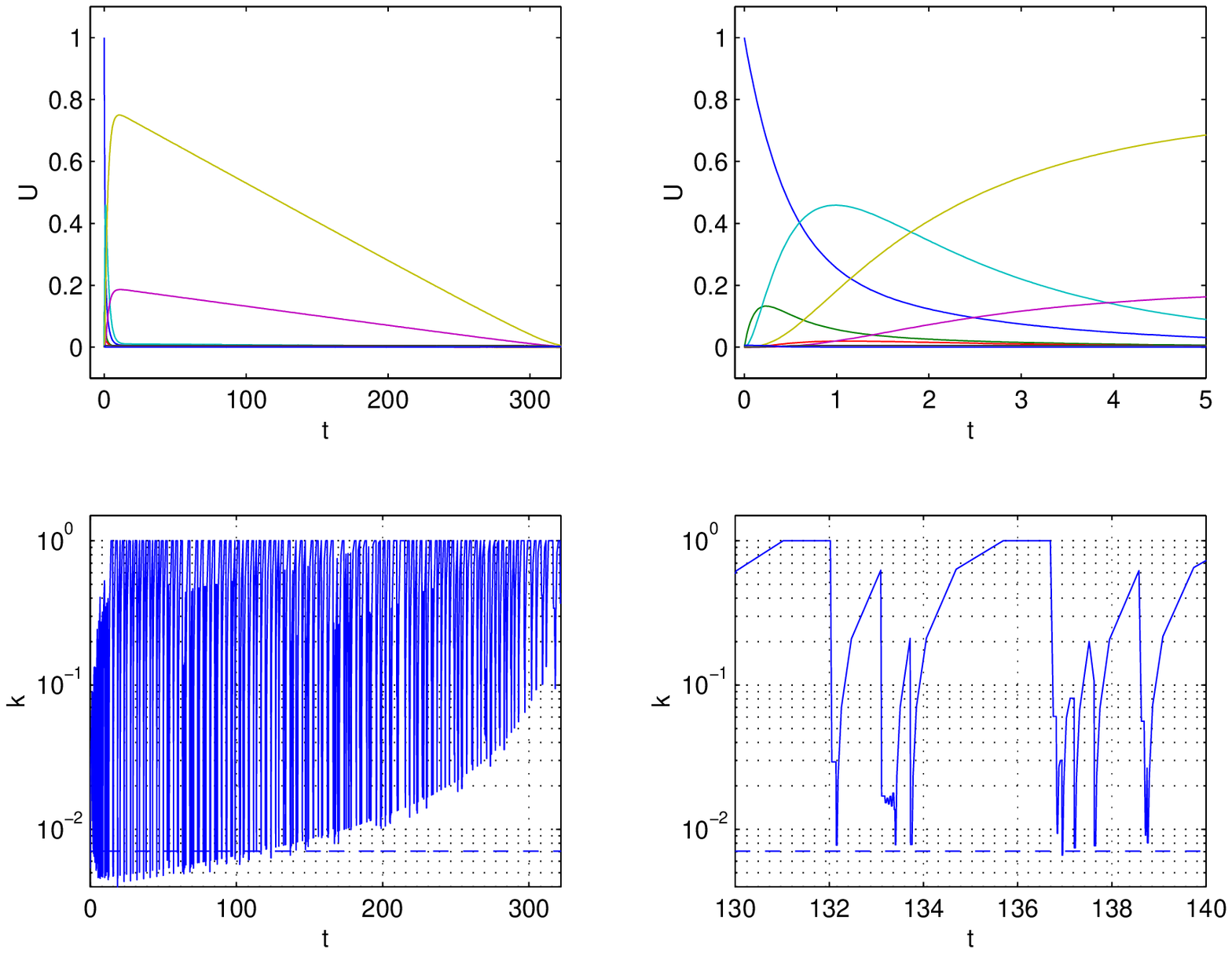}
\end{center}
\caption{Solution and time step sequence for\/ {\rm(\ref{eq:problem6})},
$\alpha/\alpha_0 \approx 1/33$.}\label{fig:problem6}
\end{figure}

The so-called HIRES problem (High Irradiance RESponse) originates from plant physiology and is
taken from the test set of ODE problems compiled by Lioen and de Swart~\cite{LioSwa98}. The problem
consists of the following eight equations:
\begin{equation}
	\label{eq:problem6}
	\left\{
	\begin{array}{rcl}
		\dot{u}_1 &=& -1.71u_1 + 0.43u_2 + 8.32u_3 + 0.0007, \\
		\dot{u}_2 &=& 1.71u_1 - 8.75u_2, \\
		\dot{u}_3 &=& -10.03u_3 + 0.43u_4 + 0.035u_5, \\
		\dot{u}_4 &=& 8.32u_2 + 1.71u_3 - 1.12u_4, \\
		\dot{u}_5 &=& -1.745u_5 + 0.43u_6 + 0.43u_7, \\
		\dot{u}_6 &=& -280.0u_6u_8 + 0.69u_4 + 1.71u_5 - 0.43u_6 + 0.69u_7, \\
		\dot{u}_7 &=& 280.0u_6u_8 - 1.81u_7, \\
		\dot{u}_8 &=& -280.0u_6u_8 + 1.81u_7,
	\end{array}
	\right.
\end{equation}
together with the initial condition $u^0 = (1.0 , 0 , 0 , 0 , 0 , 0 , 0 , 0.0057 )$.
We integrate over $[0,321.8122]$ (as specified in~\cite{LioSwa98}) and
present the solution and the time step sequence in
Figure~\ref{fig:problem6}.
The cost is now $\alpha \approx 8$, and the cost reduction factor is
$\alpha/\alpha_0 \approx 1/33$.

\begin{figure}[ht]
	\psfrag{t}{$t$}
	\psfrag{U}{$U(t)$}
	\psfrag{k}{$k(t)$}
\begin{center}
\includegraphics[width=24pc]{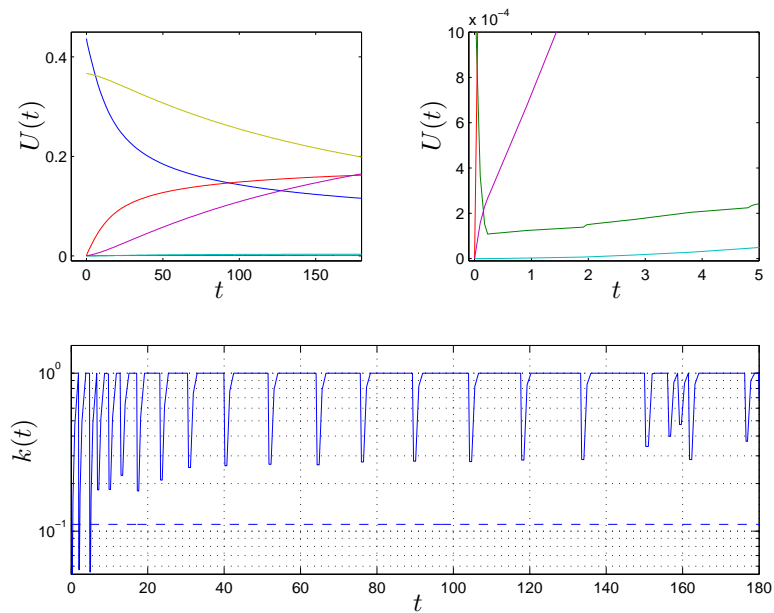}
\end{center}
\caption{Solution and time step sequence for\/ {\rm(\ref{eq:problem7})},
$\alpha/\alpha_0 \approx 1/9$.}\label{fig:problem7}
\end{figure}

\subsection{The Akzo--Nobel problem}

The next problem is a version of the ``Chemical Akzo--Nobel'' problem taken from the
ODE test set~\cite{LioSwa98}, consisting
of the following six equations:
\begin{equation}
	\label{eq:problem7}
	\left\{
	\begin{array}{rcl}
		\dot{u}_1 &=& -2r_1 + r_2 - r_3 - r_4, \\
		\dot{u}_2 &=& -0.5r_1 - r_4 - 0.5r_5 + F, \\
		\dot{u}_3 &=& r_1 - r_2 + r_3, \\
		\dot{u}_4 &=& -r_2 + r_3 - 2r_4, \\
		\dot{u}_5 &=& r_2 - r_3 + r_5, \\
		\dot{u}_6 &=& -r_5,
	\end{array}
	\right.
\end{equation}
where $F = 3.3 \cdot (0.9/737 - u_2)$ and the reaction rates are given
by $r_1 = 18.7 \cdot u_1^4 \sqrt{u_2}$, $r_2 = 0.58 \cdot u_3 u_4$,
$r_3 = 0.58/34.4 \cdot u_1 u_5$, $r_4 = 0.09 \cdot u_1 u_4^2$, and $r_5 =
0.42 \cdot u_6^2 \sqrt{u_2}$. We integrate over the interval $[0,180]$
with initial condition $u^0 = ( 0.437, 0.00123, 0, 0, 0,
0.367)$. Allowing a maximum time step of $k_{\mathrm{max}}=1$ (chosen
arbitrarily), the cost is only $\alpha \approx 2$, and the cost
reduction factor is $\alpha/\alpha_0\approx 1/9$.  The
actual gain in a specific situation is determined by the quotient
of the large time steps and the small damping time steps, as well
as the number of small damping steps that are needed. In this case, the
number of small damping steps is small, but the large time steps are
not very large compared to the small damping steps. The gain is thus
determined both by the stiff nature of the problem and the tolerance
(or the size of the maximum allowed time step).

\subsection{Van~der~Pol's equation}

\begin{figure}[t]
	\psfrag{t}{$t$}
	\psfrag{U1}{$U_1(t)$}
	\psfrag{U2}{$U_2(t)$}
	\psfrag{U3}{$U_3(t)$}
	\psfrag{k}{$k(t)$}
\begin{center}
\includegraphics[width=\figurewidth]{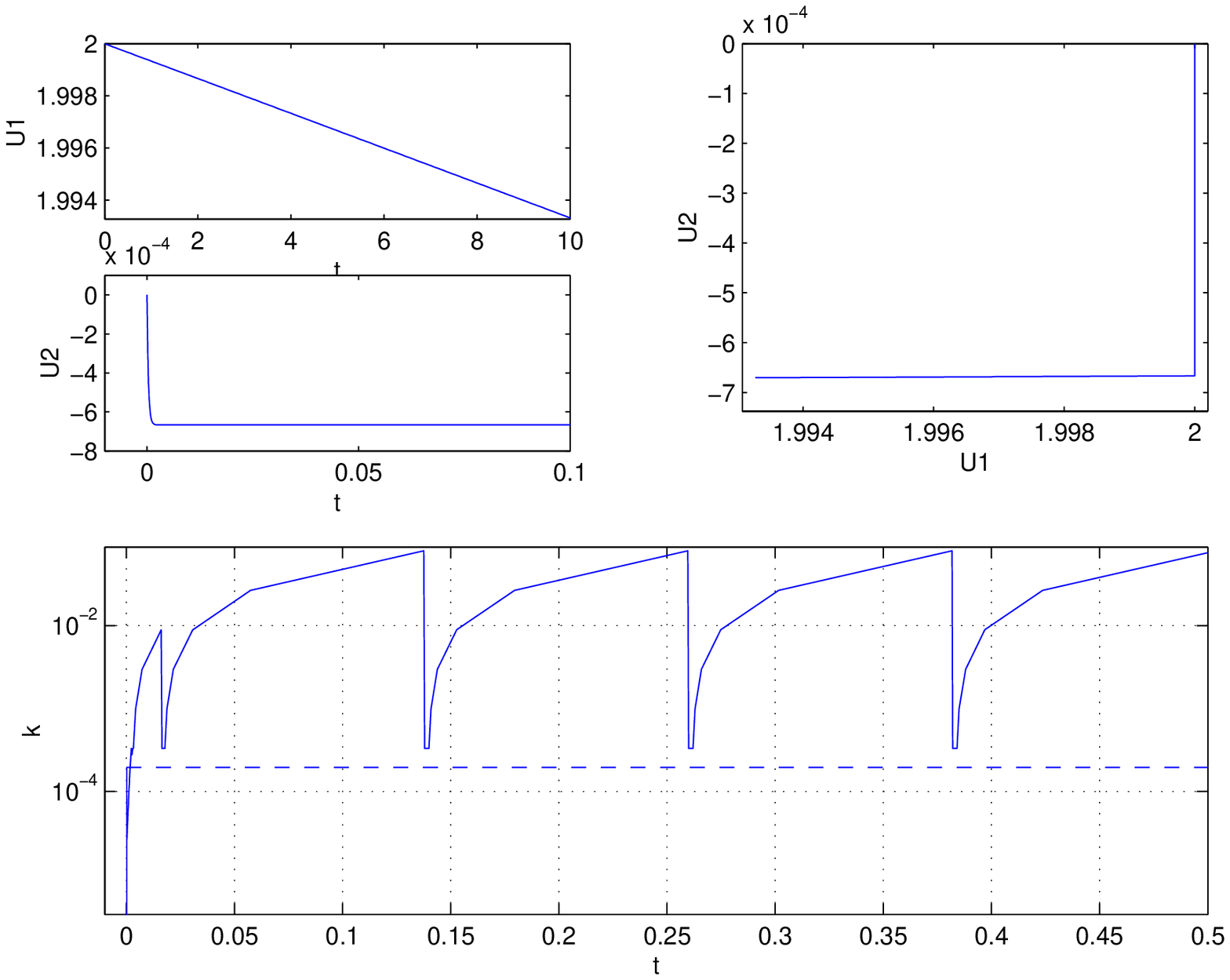}
\end{center}
\caption{Solution and time step sequence for\/ {\rm(\ref{eq:problem9})},
$\alpha/\alpha_0 \approx 1/75$.}\label{fig:problem9}
\end{figure}

A stiff problem discussed in the book by Hairer and Wanner~\cite{HaiWan91b} is
Van~der~Pol's equation,
\begin{displaymath}
	\ddot{u} + \mu (u^2-1) \dot{u} + u = 0,
\end{displaymath}
which we write as
\begin{equation}
	\label{eq:problem9}
	\left\{
	\begin{array}{rcl}
		\dot{u}_1 &=& u_2, \\
		\dot{u}_2 &=& - \mu(u_1^2 - 1)u_2 - u_1.
	\end{array}
	\right.
\end{equation}
We take $\mu=1000$ and compute the solution on the interval $[0,10]$ with initial condition
$u^0 = (2,0)$. The time step sequence behaves as desired with only a small
portion of the time interval spent on taking small damping steps (see Figure \ref{fig:problem9}).
The cost is now $\alpha \approx 140$, and the cost reduction factor
is $\alpha/\alpha_0 \approx 1/75$.

\begin{figure}[ht]
  \psfrag{t}{$t$}
  \psfrag{U}{$U(t)$}
  \psfrag{k}{$k(t)$}
\psfrag{eigenvalues}{$-\lambda_i$}
\begin{center}
\includegraphics[width=22pc]{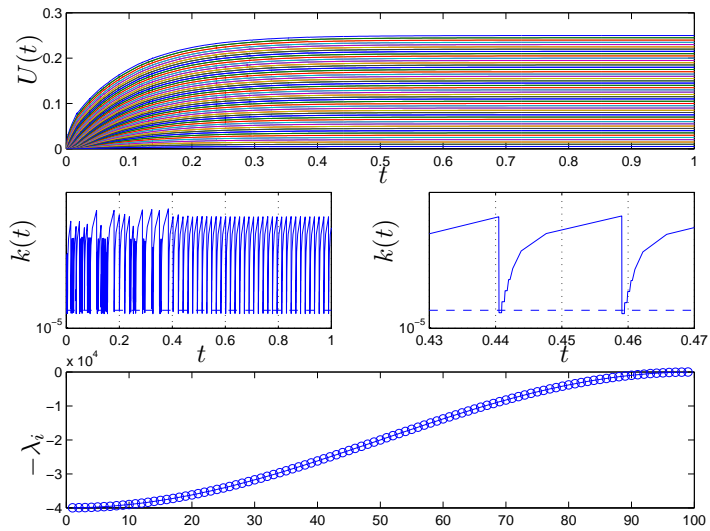}
\end{center}
\caption{Solution and time step sequence for\/ {\rm(\ref{eq:problem10})},
      $\alpha/\alpha_0 \approx 1/31$.
      Note that the time step sequence in this example is chosen slightly differently than in the previous examples, using the
      technique of dyadic damping discussed in section\/~{\rm\ref{sec:parabolic}}. This is evident upon close inspection of the
      time step sequence.}\label{fig:problem10}
\end{figure}

\begin{figure}[ht]
	\psfrag{t}{$t$}
	\psfrag{U}{$U(t)$}
	\psfrag{k}{$k(t)$}
\begin{center}
\includegraphics[width=22pc]{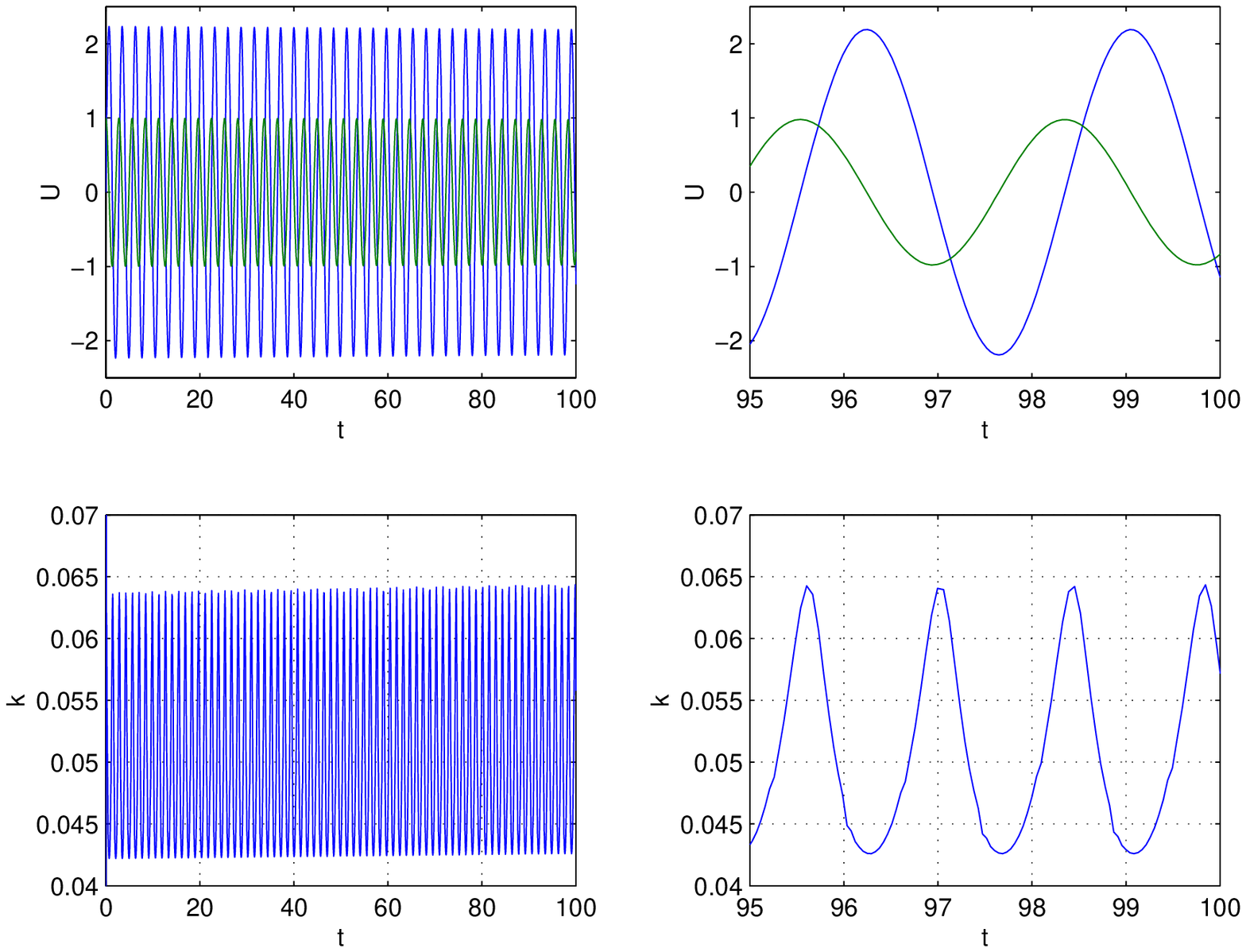}
\end{center}
\caption{Solution and time step sequence for\/ {\rm(\ref{eq:problem11})},
$\alpha/\alpha_0 = 1$.}\label{fig:problem11}
\end{figure}

\subsection{The heat equation}

A special stiff problem is the one-dimensional heat equation,
\[
\begin{split}
\dot{u}(x,t) - u''(x,t) &= f(x,t), \quad x\in(0,1), \quad t>0, \\
u(0,t) = u(1,t) &= 0, \quad t>0,\\
u(x,0)  &= 0, \quad x\in[0,1],
\end{split}
\]
where we choose $f(x,t)=f(x)$ as an approximation of the Dirac delta
function at $x=0.5$. Discretizing in space, we obtain the ODE
\begin{equation}
	\label{eq:problem10}
	\begin{split}
		\dot{u}(t) + A u(t) &= f, \\
		u(0) &= 0,
	\end{split}
\end{equation}
where $A$ is the \emph{stiffness matrix}. With a spatial resolution of
$h=0.01$, the eigenvalues of $A$ are distributed in the interval
$[0,4\cdot 10^4]$ (see Figure~\ref{fig:problem10}).

Using the dyadic damping technique described in section~\ref{sec:parabolic} for this parabolic problem,
the cost is $\alpha \approx 2000$, with a cost reduction factor of $\alpha/\alpha_0\approx 1/31$.

\subsection{A nonstiff problem}

To show that the method works equally well for nonstiff problems, we finally consider the following
simple system:
\begin{equation}
	\label{eq:problem11}
	\left\{
	\begin{array}{rcl}
		\dot{u}_1 &=& 5 u_2,\\
		\dot{u}_2 &=& -u_1.
	\end{array}
	\right.
\end{equation}
With the initial condition $u^0 = (0,1)$, the solution is $u(t) =
(\sqrt{5} \sin(\sqrt{5}t),\cos(\sqrt{5}))$.  Since this problem is
nonstiff (for reasonable tolerances), no stabilization is needed, and
so the solver works as a standard nonstiff $\mathrm{cG}(1)$ solver
with no overhead.  This is also evident from the time step sequence
(see Figure~\ref{fig:problem11}) which is chosen only to match the
size of the (continuous) residual.

\end{document}